\documentclass{amsart}
\usepackage{amsmath}
\usepackage{amssymb}
\usepackage{amsfonts}

\usepackage[foot]{amsaddr}

\makeatletter
\renewcommand{\email}[2][]{%
	\ifx\emails\@empty\relax\else{\g@addto@macro\emails{,\space}}\fi%
	\@ifnotempty{#1}{\g@addto@macro\emails{\textrm{(#1)}\space}}%
	\g@addto@macro\emails{#2}%
}
\makeatother

\begin{document}
	\title{\textbf{ON 2-ABSORBING PRIMARY HYPERIDEALS OF MULTIPLICATIVE HYPERRINGS} }

\author{Neslihan Suzen$^{1,2}$}
\author{Gursel Yesilot$^{1}$}
\address{$^{1}$ Department of Mathematics,Yildiz Technical University, 34220 Istanbul, TURKEY }
\address{$^{2}$ Department of Mathematics,University of Leicester, LE1 7RH , UK }

\email{ns433@leicester.ac.uk(N.Suzen)}
\email{gyesilot@yildiz.edu.tr (G.Yesilot)}

\maketitle

\begin{abstract}
Primary hyperideals have been introduced and studied  in multiplicative hyperrings. In this paper, we intend to study extensively  primary hyperideals of multiplicative hyperrings with absorbing zero and prove some  results regarding them. Also, we describe $ C_{u} $- ideals of multiplicative hyperrings which are  particular classes of hyperideals. In the last section, we introduce 2-absorbing primary hyperideals and investigate the properties of this notion in commutative multiplicative hyperrings.\\

\noindent {\textit{Keywords:}}multiplicative hyperring, primary hyperideal, prime hyperideal, $ C_{u} $-ideal, 2- absorbing hyperideal, 2-absorbing primary hyperideal
\end{abstract}

\section{Introduction}

The hypertructure theory was first initiated by Marty in 1934 when he defined the hypergroups \cite{marty}. Since then, algebraic hyperstructures have been investigated by many researchers with numerous applications in both pure and applied sciences. In  algebraic hyperstructures, the product of two elements is not an element but a set, while in  classical algebraic structures, the binary operation of two elements of a set  is again  an element   of the set. 	More exactly, a map $ \circ:H\times H\rightarrow P^{\ast}  (H) $ is called a \textit{hyperoperation}, where $ P^{\ast}  (H) $  is the set of all nonempty subsets  of H \cite{davvaz}.  If $ A,B\in  P^{\ast}  (H) $ and $ x\in H $, then we define 
$$ A\circ B=\bigcup_{a\in A, b\in B} a\circ b,\quad  A\circ x = A\circ\{x\}
$$ A \textit{semihypergroup} $ (H,\circ) $ is a nonempty set with  the associative hyperoperation; i.e   $ a \circ (b\circ c)=(a \circ b) \circ c $ for all $ a,b,c \in H $. A semihypergroup $ H $ is called a \textit{hypergroup} if for every $ a\in H $, $ a \circ H=H=H \circ a $, which is \textit{quasihypergroup}.

Similar to hypergroups,  hyperrings are algebraic structures more general than  rings, subsitutiting both or only one of the binary operations of addition and multiplication by hyperoperations. The hyperrings were  introduced  by Krasner  \cite{krasner}. \textit{Krasner hyperrings} are a generalization of  classical rings in which the multiplication is a binary operation while the addition is a hyperoperation. Since then, this concept has been studied by many researchers. The another type of hyperrings called  Multiplicative hyperring was introduced and studied by Rota in 1982 \cite{rota}, which was subsequently investigated by many authors [4,5,8,10]. A \textit{multiplicative hyperring} is a hyperstructure $ (R,+,\cdot) $, where $ (R,+) $ is  an additive commutative group  and $ (R,\cdot)  $ is a semihypergroup which satiesfies the axioms   i)  $ x \cdot (y+z)\subseteq x \cdot y+x \cdot z  $ and $ (y+z) \cdot x \subseteq y \cdot x+z \cdot x$   and 
ii)    $ x \cdot (-y)=(-x) \cdot  y=-(x \cdot y) $ for all   $x,y,z\in R $. 
If in (i) we have equalities instead of inclusion, then we say that the multiplicative hyperring is \textit{strongly distributive}. 
 A hyperring  $ R $ is called \textit{proper hyperring} if it is not a ring. Also, $ R $ is said to be \textit{commutative} if $ R $ is commutative with respect to operation $ + $ and hyperoperation $ \cdot $ .  
 
 Recall from \cite{ameri} that a nonempty subset $ I $  of a commutative  hyperring $ R $ is said to be a \textit{hyperideal} of $ R $ if  $x-y\in I $  and $ r\cdot x\subseteq I $ for any  $ x,y\in I  $ and $ r\in R $. If $ A $ and $ B $ are hyperideals of $ R $, then $ A+B=\{a+b: a\in A,b\in B\}  $ and $ AB=\bigcup\{ \sum_{i=1}^{n}a_{i}\cdot b_{i}:a_{i}\in A, b_{i}\in B  \quad and \quad n\in \mathbb{N}\} $, which are hyperideals of $ R $ [4,5]. This was defined in \cite{dasgupta} that the  \textit{principal hyperideal } of  $ R $ generated by an element $ a $ is given by $ <a>=\{pa:p\in \mathbb{Z}\}+\{\sum_{i=1}^{n}x_{i}+\sum_{j=1}^{m}y_{j}+\sum_{k=1}^{l}z_{k}:\forall i,j,k,\exists r_{i},s_{j},u_{k}\in R \quad such\quad that\quad x_{i}\in r_{i}\circ a, y_{j}\in a\circ s_{j}, z_{k}\in t_{k}\circ a\circ u_{k}\}  $. We define a set  $A=\{r\in R:r^{n}\subseteq I\quad for\quad some\quad n\in \mathbb{N}\} $, is called as the \textit{radical of $ I $} and denoted by $ \sqrt{I} $ or $ Rad(I) $. A proper hyperideal $ P $ of $ R $ is called a \textit{prime hyperideal} if for any $ x,y\in R $,  $ x\cdot y\subseteq  P$ then  $ x\in P $ or $ y\in P $ \cite{procesi}.  Recall that a proper ideal $ Q $ of $ R $ is said to be a \textit{primary hyperideal} of $ R $ if  $ x\cdot y\subseteq Q  $, then $ x\in Q $ or $ y\in \sqrt{Q} $ for any $ x,y\in R $ \cite{dasgupta}. Every prime hyperideal of a commutative multiplicative hyperring is clearly a primary hyperideal. However, the converse is not true in general[Example 3.5, 5]. 

 The concept of $ C $-ideal of a multiplicative hyperring was defined by Dasgupta in 2012 \cite{dasgupta}, which is a particular class of hyperideals. An hyperideal $ I $ of a multiplicative hyperring $ R  $ is said to be a \textit{$ C $-ideal } if  whenever $ C $ is the class of all finite products of elements of $ R $, i.e; $  C=\{r_{1}\cdot r_{2}\cdot...\cdot r_{n}:r_{i}\in R,n\in\mathbb{N}\}\subseteq P^{*}(R	) $ and $ A\cap I\neq\varnothing $, then $ A\subseteq I $  for any  $  A\in C  $. Various generalizations of prime hyperideals, primary hyperideals and also $ C $-ideals of a multiplicative  hyperring were studied in his paper. Note that if $ Q $ is a primary $ C $-ideal of a hyperring $ R $, then $ Rad(I) $  is a prime hyperideal [Proposition 3.6, 5].
 
 Let  $ (R_{1},+, \circ _{1} ) $ and  $ (R_{2},+,\circ_{2}   ) $ be  hyperrings and  $ f :R_{1}\longrightarrow R_{2} $ be a map. Then $ f $ is called a \textit{homomorphism( or good homomorphism)} if for all  $ x,y\in R_{1} $, $ f(x+y)=f(x)+_{2}f(y) $  and $ f(x\circ y)\subseteq f(x)\circ_{2}f(y) $ (or $ f(x\circ y)= f(x)\circ_{2}f(y) $). Also, let $ f $ is a  good homomorphism. \textit{The kernel of f} is the inverse image of $ <0> $, the hyperideal generated by the zero in $ R_{2} $  and it is denoted by $ Kerf $ \cite{davvaz}. Clearly, we have $ f(<0>)\subseteq <0> $, which follows   $ <0>\subseteq Kerf $. Note that $ a-b\in Kerf  $ does not imply $ f(a)=f(b) $ in hyperstructures since the zero hyperideal can contain more than 0, where the \textit{zero hyperideal}   is the hyperideal generated by the additive identity 0;  i.e $ <0>=\{\sum_{i}x_{i}+\sum_{j}y_{j}+\sum_{k}z_{k}:each\quad sum \quad is\quad finite\quad and\quad for\quad each \quad i,j,k\quad there\quad exist\quad r_{i},s_{j},t_{k},u_{k}\in R\\  such \quad that\quad x_{i}\in r_{i}\cdot0,y_{i}\in0\cdot s_{j},z_{k}\in t_{k}\cdot0\cdot u_{k} \} $ \cite{davvaz}.
 
 Recall that we can define \textit{quotient multiplicative hyperrings} similar to quotient rings in classical algebra \cite{davvaz}. Let $ R $ be a multiplicative hyperring and $ I $ be a hyperideal of $ R $. We consider  the usual addition of cosets  and multiplication $ (a+I)\ast(b+I)=\{c+I:c\in a.b\} $ on the set  $ R/I=\{a+I: a\in R\} $ of all cosets of $ I $. Then $ (R/I,+,\ast) $ is a multiplicative  hyperring. Note that if $ I $ and $ K $ are hyperideals of a multiplicative hyperring $ R $ such that $ K\subseteq I $, then $ R/K $ is a multiplicative hyperring and $ I/K $ is also a hyperideal of   $ R/K $. Some results regarding quotient hyperrings can be seen in \cite{davvaz} .
 
 In this paper, we obtain some results and examples about prime and primary hyperideals of multiplicative  hyperrings.  Also, we introduce and study \textit{minimal prime hyperideal} and \textit{$ C $-union ideal} in commutative multiplicative hyperrings, which is a particular class of hyperideals.  It is shown that $ C $-ideal and $ C $-union ideal are different concepts. We prove that every $ C $-union ideal is a $ C $-ideal of hyperring $  R $, but the converse is not true. We study  homomorphisms of  hyperring regarding $ C $-union ideals. It is shown that the surjective homomorphic image of a prime $ C $-union ideal (primary $ C $-union ideal) hyperideal  is also a prime(primary) hyperideal  on the  condition that it contains the  kernel of homomorphism. In the other section, we define \textit{2-absorbing hyperideals} and  \textit{2- absorbing primary  hyperideals} of a multiplicative hyperring, which were studied by Badawi, Tekir  and Yetkin  in ordinary algebra [2,3]. It is shown that every 2-absorbing hyperideal is a 2- absorbing primary hyperideal, but the converse is not true. Also, we prove that every primary hyperideal is 2- absorbing primary hyperideal. It is shown that  a 2- absorbing primary hyperideal need not to be a primary hyperideal. We investigate   results and examples that show some differences from same notions in ordinary algebra. 
 
 Throughout this paper, we assume that all hyperrings are proper commutative multiplictive hyperrings with absorbing zero; i.e    there exists $ 0\in R $  such that x=0+x and $ 0\in x\cdot0=0\cdot x $ for all $ x\in R $.

\section{Properties of Hyperideals}
\textbf{Lemma 2.1.} Let $ I $ be a prime hyperideal of hyperring $ R $ and $ J $ be a subset of $ R $. For any $ a\in R $  , $ aJ\subseteq I $ and $ a\notin I $ imply that $ J\subseteq I $.\\
\\
\textbf{Proof.} Let  $ aJ\subseteq I $ and $ a\notin I $ for any $ a\in R $. Hence, we have   $ aJ= \cup aj_{i}  \subseteq I $ for all $ j_{i}\in J $. Then,  $ aj_{i} \subseteq I $ for all $ j_{i}\in J $. Since   $ I $ prime hyperideal and $ a\notin I $, we conclude that  $ j_{i}\in I $ for all $ j_{i}\in J $. Thus  $ J\subseteq I $.$ \Box $\\ 
\\
\textbf{Lemma 2.2.} Let $ I $ be a primary hyperideal of $ R $ and  $ J $ be a subset of $ R $. For any $ a\in R $, $ aJ\subseteq I $ and $ a\notin I $ implies that $ J\subseteq \sqrt{I} $.(or $ aJ\subseteq I $ and $ J\nsubseteq I $ imply that $ a\in \sqrt{I} $). \\
\\
\textbf{Proof.} Let $ aJ\subseteq I $ and $ a\notin I $ for any $ a\in R $. Then we have   $  aJ=\cup aj_{i}  \subseteq I $ for all $ j_{i}\in J $. Hence   $ aj_{i} \subseteq I $ for all $ j_{i}\in J $. Since   $ I $ primary hyperideal and $ a\notin I $, we conclude that    $ j_{i}\in \sqrt{I} $ for all  $ j_{i}\in J $. Thus $ J\subseteq \sqrt{I} $. The proof of the other argument is similar.$ \Box $\\
\\
\textbf{Proposition 2.3.} Let $ I $ be a prime hyperideal of   $ R $ and $ A $, $ B $ be subsets of $ R $. If  $ AB\subseteq I $, then $ A\subseteq I $ or $ B\subseteq I $.\\
\\
\textbf{Proof.} Suppose that $ AB\subseteq I $,  $ A\nsubseteq I $ and $ B\nsubseteq I $. Since  $ AB=\bigcup a_{i}b_{i}\subseteq I $ , we have $ a_{i} b_{i}\subseteq I $ for all $ a_{i}\in A, b_{i}\in B $.Since  $ A\nsubseteq I $ and $ B\nsubseteq I $, then there exist $ x,y\notin I $ for some  $ x\in A $, $ y\in B $. This implies    $ xy\subseteq AB\subseteq I $ .Since $ x,y\notin I $ and $ I $ is a  prime hyperideal, then  $ xy\nsubseteq I $, a contradiction.Thus $ A\subseteq I $ or $ B\subseteq I $.$ \Box $\\
\\
\textbf{Definition 2.4.}  Let $ I$ be a hyperideal of hyperring $ R  $ and    $ P $ be a prime  hyperideal such that $ I\subseteq P $ . If there is no prime hyperideal $  P^{'} $ such that    $ I\subseteq P\subseteq P^{'} $, then $ P $ is called  \textit{minimal prime hyperideal of $ I $}. The set of all minimal prime hyperideals of $ I $ is denoted by $  Min_{h}(I) $.\\
\\
\textbf{Proposition 2.5.} If  $  P $ is a prime hyperideal of  $ R $ then $ Min_{h}(P)=\{P\} $.\\
\\
\textbf{Proof.} Since $ P $ is a prime hyperideal, the proof is clear .$ \Box $\\
\\
\textbf{ Example 2.6.} Let $ (\mathbb{Z},+,\cdot) $ be the ring of integers. For all $ x,y\in\mathbb{Z} $; we define the hyperoperation $ x\circ y=\{2xy,3xy\} $.   Then  $ (\mathbb{Z},+,\circ) $ is a multiplicative hyperring.     The set  $ 12\mathbb{Z}=\{12n:n\in\mathbb{Z}\} $ is a hyperideal that is not prime hyperideal of $ R $. Moreover, $ 2\mathbb{Z}  $ and $ 3\mathbb{Z}   $ are minimal prime hyperideals of $ 12\mathbb{Z} $.\\
\\
\textbf{Proposition 2.7.} Let $ I_{1},I_{2},... ,I_{n} $ be some hyperideals of a hyperring  $ R $. Then the following statements hold.
\begin{enumerate}
	\item $ Rad(I_{1}I_{2}...I_{n})=Rad(\bigcap^{n} _{i=1}I_{i})=\bigcap_{i=1}^{n}Rad(I_{i}) $ 
	\item $ Rad(I)\subseteq Rad(Rad(I)) $ 
\end{enumerate} 
\textbf{Proof.} 
\begin{enumerate} 
	\item Since $ I_{1}I_{2}...I_{n}\subseteq \bigcap^{n} _{i=1}I_{i} $, $ Rad(I_{1}I_{2}...I_{n})\subseteq Rad(\bigcap^{n} _{i=1}I_{i}) $. Now let $ x\in Rad(\bigcap^{n} _{i=1}I_{i}) $. Then there exists $ m\in \mathbb{N} $ such that  $ x^{m}\subseteq \bigcap^{n} _{i=1}I_{i} $. Thus,  $ x^{m}\subseteq I_{i} $ for all  $ i=1,...,n  $. Hence for all  $ i $,  $ x\in Rad(I_{i}) $ and then $ x\in \bigcap_{i=1}^{n}Rad(I_{i}) $. It follows $  Rad(\bigcap^{n} _{i=1}I_{i})\subseteq\bigcap_{i=1}^{n}Rad(I_{i}) $. Conversely, suppose that $ x\in \bigcap_{i=1}^{n}Rad(I_{i})$. Then for all $ i $, $  x\in Rad(I_{i})$. Thus there exists $ m_{i}\in \mathbb{N }$ such that  $ x^{m_{i}}\subseteq I_{i} $. Since $ x^{m_{1}}.x^{m_{2}}...x^{m_{n}}\subseteq I_{1}I_{2}...I_{n} $, we conclude that $ x^{\sum_{i=1}^{n}}\subseteq I_{1}I_{2}...I_{n} $ and so  $  x\in Rad (I_{1}I_{2}...I_{n}) $. Hence $ \bigcap_{i=1}^{n}Rad(I_{i})\subseteq Rad(I_{1}I_{2}...I_{n}) $. Now let  $ x\in Rad(\bigcap^{n} _{i=1}I_{i}) $, then there exists $ m\in \mathbb{N} $ such that   $ x^{m}\subseteq \bigcap^{n} _{i=1}I_{i} $. Thus  $ x^{m}\subseteq I_{i} $ $ x^{mn}\subseteq I_{1}I_{2}...I_{n} $ for all $ i $ and so  $ x\in Rad(I_{1}I_{2}...I_{n})  $.
	\item Since $ I\subseteq Rad(I) $ for all hyperideal $ I $, the proof is clear . $ \Box $
\end{enumerate}
If  $ I $ is a  $ C $-ideal, in (2) we have equality instead of inclusion \cite{dasgupta}.\\
\\
\textbf{Proposition 2.8.} Let $ f:R\longrightarrow S $ be an onto good  homomorphism of hyperrings. If  $ I $ is a hyperideal of  $ R $, then $ f(\sqrt{I}) \subseteq \sqrt{f(I)}$ .\\
\\
\textbf{Proof.} Let  $ x\in f(\sqrt{I}) $ for any $ x\in R $. Since  $ f $ is onto, then there exists $ y\in \sqrt{I} $ such that $ x=f(y) $. It implies that $ y^{n}\subseteq I $ for some  $ n\geq1 $. Since  $ f $ is a good homomorphism, we conclude that  $ x^{n}=f(y)^{n}=f(y^{n})\subseteq f(I) $. Thus $ x\in \sqrt{f(I)} $.$ \Box $\\
\\
\textbf{Proposition 2.9.} Let $ f:R\longrightarrow S $ be good homomorphism of hyperrings and let $ J $ be a hyperideal of $ S $. Then  $ \sqrt{f^{-1}(J)} =f^{-1}(\sqrt{J})$.\\
\\
\textbf{Proof.} Let $ r\in   f^{-1}(\sqrt{J}) $. Then   $ f(r)\in \sqrt{J} $. Since  $ f  $ is a good homomorphism, for some $ n\geq  1$, we have  $ (f(r))^{n}=f(r^{n})\subseteq J $. Hence $ r^{n}\subseteq f^{-1}(J) $ and so $ r\in \sqrt{f^{-1}(J)} $. The converse can be shown similarly .$ \Box $\\
\\
\textbf{Definition 2.10.} Let $ C $ be the  class of all finite hyperproducts of  elements of multiplicative hyperring $ R $ ,i.e, $ C=\{r_{1}\cdot r_{2}\cdot...\cdot r_{n}:r_{i}\in R,i=1,...,n\} $. If for any  $ A_{j}\in C $,  $ (\cup A_{j})\cap I \neq\varnothing $ implies that  $ \cup A_{j}\subseteq I $, then $ I $ is said to be a \textit{$ C $-union ideal} of $ R $ and it denotes by  $ C_{u} $-ideal.\\
\\
\textbf{Example 2.11.} Let $ (\mathbb{Z},+,\cdot)  $ be the ring of integers. We define the hyperoperation   $ x\circ y=\{2xy,4xy\} $ for all  $ x,y\in\mathbb{Z} $, then  $ (\mathbb{Z},+,\circ) $ is a multiplicative hyperring. In $ (\mathbb{Z},+,\circ) $,  since  $ x\circ y=\{2xy,4xy\} $ for any $ x,y\in R $, then all finite products  of elements are subsets of the hyperideal $ 2\mathbb{Z}=\{2n:n\in \mathbb{Z}\} $. Since for all finite products $ A_{j} $,  $ (\cup A_{j}) \cap I\neq\varnothing $ and  $ \cup A_{j}\subseteq2\mathbb{Z} $, then the hyperideal $ 2\mathbb{Z} $ is a $ C_{u} $-ideal.\\
\\
Note that every $ C_{u} $-ideal of $ R $ is a  $ C $-ideal. In fact, let $ I $ be a  $ C_{u} $-ideal and for some $ A_{j}\in C $ ,  $ A_{j}\cap I\neq\varnothing $.  Hence, for any union $ \cup A_{j} $ of elements of $ C $ such that $ A_{j}\cap I\neq\varnothing $,  we have  $ (\cup A_{j})\cap I	\neq\varnothing  $. Since   $ I $ is a  $ C_{u} $-ideal, then $ \cup A_{j} \subseteq I $ and so  $ A_{j}\subseteq I $.\\
\\
It is easy to see that $ C $-ideal and $ C_{u} $-ideal are different concepts. The following is an example of a hyperideal where it is a $ C $-ideal, but it is not a $ C_{u} $-ideal.\\
\\
\textbf{Example 2.12} Consider the hyperring $ (\mathbb{Z},+,\circ) $ in  Example 2.6. The hyperideal $ 5\mathbb{Z}=\{5n:n\in \mathbb{Z}\} $  of $ \mathbb{Z} $ is a  $ C $-ideal, but  it is not a $  C_{u} $-ideal. Clearly $ (1\circ 1 \cup5\circ 1)\cap5\mathbb{Z}=\{10,15\}\neq\varnothing $ ,but   $ 1\circ 1 \cup5\circ 1=\{2,3,10,15\}\nsubseteq 5\mathbb{Z} $.\\
\\
Recall that if $ Q $ is a hyperideal of  hyperring $ S $ and $ f:R\longrightarrow S $ is a good homomorphism, then $ f^{-1}(Q) $ is always a hyperideal of hyperring $ R $. However, if $ I  $ is a hyperideal of the hyperring $ R $ and $ f:R\longrightarrow S $ is a good homomorphism, then $ f(I) $ need not to be  a hyperideal of $ S  $. \\
\\ 
\textbf{Theorem 2.13.} Let $ f:R\longrightarrow S $ be a homomorphism of hyperrings. Then the following statements hold.
\begin{enumerate}
	\item Let $ I $ is a $ C_{u} $- ideal of  hyperring  $  R $ and  $ f:R\longrightarrow S $ is an onto good homorphism such that  $ Kerf\subseteq I $. If  $ I $ is a prime hyperideal  of $ R $, then $ f(I) $ is a prime hyperideal of  hyperring  $ S $.
		\item If $ f:R \longrightarrow S $  is a good homomorphism and $ J$ is a prime hyperideal of  hyperring $ S $, then  $ f^{-1}(J ) $ is a prime hyperideal of  hyperring $ R $. 
\end{enumerate}
\textbf{Proof.}\\
\textbf{1.} It is easy to see that $ f(I ) $ is a hyperideal of $ S $. Now, let $ s_{1}.s_{2}\subseteq f(I) $ for any $ s_{1},s_{2}\in S $. Since $ f $ is onto homomorphism, then there exists $ r_{1}, r_{2}\in R $ such that    $ s_{1}=f(r_{1}) $ and $ s_{2}=f(r_{2}) $.  Thus, $f(r_{1}.r_{2})= f(r_{1}).f(r_{2})=s_{1}.s_{2}\subseteq f(I) $  and so $ 0\in f(I)-f(r_{1}.r_{2})=f(I-r_{1}.r_{2})=\{f(u):u\in I-r_{1}.r_{2}\}  $. Then, there exists   $ v\in I-r_{1}.r_{2} $ such that  $ f(v)=0\in<0> $ and hence  $ v\in Kerf $. Since   $ (I-r_{1}.r_{2})\cap I\neq\varnothing $ and $ I $ is a  $ C_{u} $- ideal, we conclude  $ I-r_{1}.r_{2}\subseteq I $. Thus, $ r_{1}.r_{2}\subseteq I $. Since $ I $ is a prime hyperideal, it implies that  
$ r_{1} \in I $ or $ r_{2}\in I $. and so $ s_{1}\in f(I) $ or $ s_{2}\in f(I) $. Therefore, $ f(I) $ is a prime hyperideal of $ S $.\\
\textbf{2.} Let  $r_{1},r_{2} \in R $ such that $ r_{1}r_{2} \subseteq f^{-1}(J) $. Since $ f(r_{1}r_{2})=f(r_{1})f(r_{2})\subseteq J $ and $ J $ is a prime hyperideal of $ S $, we have  $ f(r_{1})\subseteq J $ or $ f(r_{2})\subseteq {J} $. Hence we conclude that
$ r_{1} \in f^{-1}(J) $ or $ r_{2} \in f^{-1}(J) $. Thus $ f^{-1}(J) $ is a prime hyperideal of $ R $.$ \Box $\\
\\
\textbf{Theorem 2.14.} Let  $ f:R\longrightarrow S $ be a homomorphism of hyperrings. Then the following statements hold.
\begin{enumerate}
	\item Let $ f $ is an onto good homomorphism and   $ I $ be a $ C_{u} $- ideal of  $  R $ such that  $ Kerf\subseteq I $. If   $ I $ is a primary hyperideal, then   $ f(I)$ is a primary hyperideal of $ S $.
	\item If $ f:R \longrightarrow S $ is a good homomorphism of hyperrings and  $ J$ is a primary hyperideal of  $ S $, then  $ f^{-1}(J )$  is a primary hyperideal of  $R $.
\end{enumerate}
\textbf{Proof.}\\
\textbf{1.} Let $ s_{1}s_{2}\subseteq f(I ) $ for any  $ s_{1},s_{2}\in S $. Since  $ f $ is onto homomorphism, then there exists $ r_{1},r_{2}\in R $ such that  $ s_{1}=f(r_{1}) $ and $ s_{2}=f(r_{2}) $. Hence, $  f(r_{1}r_{2})= f(r_{1})f(r_{2})=s_{1}s_{2}\subseteq f(I)  $ implies  $ f(r_{1}r_{2})\subseteq f(I) $ and so $ 0\in  f(I)-f(r_{1}r_{2})=f(I-r_{1}r_{2})=\{f(u):u\in I-r_{1}r_{2}\} $. Thus    $ f(v)=0\in <0> $ for some $ v\in I-r_{1}r_{2} $. It follows that  $ v\in Kerf\subseteq I $. Since $ I\cap (I-r_{1}r_{2})\neq\varnothing $ and $ I $ is a $ C_{u} $-ideal, we have   $ I-r_{1}r_{2}\subseteq I $. Since $ r_{1}r_{2}\subseteq I $ and $ I $ is a primary hyperideal of $ R $, we conclude that $ r_{1}\subseteq \sqrt{I}  $ or $ r_{2} \subseteq \sqrt{I} $. Then  $ s_{1}\subseteq f(\sqrt{I})\subseteq\sqrt{f(I)} $ or $ s_{2}\subseteq f(\sqrt{I })\subseteq\sqrt{f(I)}   $ by Proposition 2.8.\\
\textbf{2.} Let   $ xy\subseteq f^{-1}(J) $ for any  $ x,y\in R $. Since  $ f(xy)=f(x)f(y)\subseteq J $ and $ J $ is a primary hyperideal of $ S $, we get $ f(x)\subseteq J $ or $ f(y)\subseteq \sqrt{J} $. Hence,  
$f(x)\subseteq J $ or $ f(y)\subseteq \sqrt{J} $ and so  $ x\subseteq f^{-1}(J) $ or $ y\subseteq f^{-1}(\sqrt{J}) $. Since    $\sqrt{f^{-1}(J)}=f^{-1}(\sqrt{J})  $,  we conclude that  $f^{-1}(J) $ is a primary hyperideal of $ R $  by Proposition 2.9.$\Box $\\	
\\
\section{2-Absorbing Primary Hyperideals of Hyperrings}
\textbf{Definition 3.1.}  Let   $ I $ be a proper hyperideal of a hyperring $ R $.  The hyperideal $ I $ is called \textit {2-absorbing hyperideal} of $ R $  if   $ a\cdot b\cdot c\subseteq I $, then $ a\cdot b\subseteq I $ or $ b\cdot c\subseteq I $ or $ a\cdot c\subseteq I $ for any $ a,b,c\in R $.\\
\\ 
\textbf{Definition 3.2.} Let   $ I $ be a proper hyperideal of a hyperring  $ R $. The hyperideal $ I $ is called \textit{2-absorbing primary hyperideal} of $ R $  if  $ a\cdot b\cdot c\subseteq I $, then $ a\cdot b\subseteq I $ or $ b\cdot c\subseteq \sqrt{I} $ or $ a\cdot c\subseteq\sqrt{I} $ for any  $ a,b,c\in R $. \\
\\
It is clear that every 2-absorbing hyperideal is a 2-absorbing primary hyperideal. The converse is not true, as is  shown in the following example.\\
\\
\textbf{Example 3.3.}\\
\textbf{(1) } Let $ (\mathbb{Z},+,\cdot)$ be the ring of integers. For all $ x,y\in\mathbb{Z} $; we define the hyperoperation $ x\circ y=\{2xy,3xy\} $.   Then  $ (\mathbb{Z},+,\circ) $ is a multiplicative hyperring. The subset  $ 12\mathbb{Z}=\{12n:n\in\mathbb{Z}\} $ is a 2-absorbing primary hyperideal that is not  2-absorbing hyperideal of $\mathbb{Z}  $. \\
\\
\textbf{(2)} Consider  the ring of integers $ \mathbb{Z} $. For all  $  x,y\in\mathbb{Z} $; we define the hyperoperation $ x\circ y=\{2xy,4xy\} $. Then $ (\mathbb{Z},+,\circ) $ is a multiplicative hyperring. The hyperideal  $ 120\mathbb{Z}=\{120n:n\in\mathbb{Z}\} $ is a 2-absorbing primary hyperideal,  but it is not a   2-absorbing hyperideal. Also, the hyperideal $ 15\mathbb{Z}=\{15n:n\in\mathbb{Z}\} $ is a  2-absorbing  hyperideal.\\
\\
\textbf{(3)} Consider the ring $ (\mathbb {Z}_6,\oplus,\odot)  $ that $ \bar{a}\oplus\bar{b} $ and $ \bar{a}\odot\bar{b} $ are remainder of  $ \frac{a+b}{6} $ and $ \frac{a\cdot b}{6} $  which + and $ \cdot $ are ordinary addition and multiplication for all $ \bar{a},\bar{b}\in\mathbb{Z}_{6} $. For all $  \bar{a} ,\bar{b}  \in \mathbb{Z}_6 $, we define  the hyperoperation  $ \bar{a} \star \bar{b}=\{\overline{ab} ,\overline{2ab},\overline{3ab},\overline{4ab},\overline{5ab}\} $. Then  $ (\mathbb {Z}_6,\oplus,\star)  $ is a commutative multiplicative hyperring. The hyperideal  $  \{ \bar{0}\} $ of $ (\mathbb {Z}_6,\oplus,\star)  $ is  2-absorbing hyperideal.\\
\\ 
Note that every primary hyperideal is 2-aborbing primary hyperideal. In fact; let  $  I$ be a primary hyperideal of $R $. Suppose that     $ abc\subseteq I $ and  $ ab\nsubseteq I $ for any $ a,b,c\in R $. Since $ I $ is a primary hyperideal, then  $ c\subseteq \sqrt{I} $  by lemma 2.2. Hence there exists $ n>0  $ such that $ c^{n}\subseteq I $. Since $ I $ is a hyperideal, we have  $ a^{n}c^{n}\subseteq I $ and $ b^{n}c^{n}\subseteq I $. Thus $ ac\subseteq\sqrt{I} $ and $ bc\subseteq\sqrt{I} $ and so  $ I $ is  2-absorbing primary hyperideal.\\
\\
The following example shows that a 2-aborbing primary hyperideal need not to be a primary hyperideal.\\
\\
\textbf{Example 3.4.}\\
\textbf{(1)} Consider the hyperring $ (\mathbb{Z},+,\circ) $ in Example 3.3\textbf{(1)}. The hyperideal  $ 12\mathbb{Z}=\{12n:n\in\mathbb{Z}\}  $ is a 2-absorbing primary hyperideal of $ \mathbb{Z} $ , but $ 12\mathbb{Z} $ is not a primary hyperideal of $ \mathbb{Z} $. Clearly,  $4\circ3\subseteq12\mathbb{Z} $ and $ 4\notin12\mathbb{Z} $, but  for all $ n>0 $ we have  $ 3^{n}\nsubseteq12\mathbb{Z}  $ ; and $4\circ3\subseteq12\mathbb{Z}$ and $ 3\notin12\mathbb{Z} $, but for all $ n>0 $  we have $ 4^{n}\nsubseteq12\mathbb{Z}  $. Thus  $ 12\mathbb{Z} $ is not primary hyperideal.\\
\\
\textbf{(2)}Consider the hyperring $ (\mathbb{Z},+,\circ) $ in  Example 3.3\textbf{(3)}. Since the hyperideal  $  \{ \bar{0}\} $ is 2-absorbing hyperideal, then it is 2- absorbing primary hyperideal. However, the hyperideal  $  \{ \bar{0}\} $ is not a primary hyperideal. In fact,  $  \bar{3}\star \bar{2}=\{\bar{0}\}  $, but $ \bar{3}\neq\bar{0} $ and $ \bar{2}\neq\bar{0} $. Also, we have   $ \bar{2}^{n}=\{\bar{0},\bar{2},\bar{4}\}\nsubseteq\{\bar{0}\} $  and $ \bar{3}^{n}=\{\bar{0},\bar{3}\}\nsubseteq\{\bar{0}\}  $ for all  $ n>0 $. Thus,  $ \{\bar{0}\} $ is not a primary hyperideal. \\
\\
\textbf{Theorem 3.5.} Let $ I $ be a hyperideal of the hyperring $ R $. If $ \sqrt{I} $ is a prime hyperideal, then  $ I $ is a 2-absorbing primary hyperideal of  $ R $.\\
\\
\textbf{Proof.} Suppose that  $ abc\subseteq I $ and  $ ab\nsubseteq I $ for any  $ a,b,c\in R $. Since $ (ac)(bc)=abc^{2}\subseteq I\subseteq\sqrt{I} $ and  $ \sqrt{I} $ is a prime hyperideal, we have     $ ac\subseteq\sqrt{I} $ or $ bc\subseteq\sqrt{I} $ by Proposition 2.3. Hence, $ I $ is a  2-absorbing primary hyperideal of $ R $. $\Box$\\
\\
The converse of  Theorem 3.5 is not true in general .\\
\\
\textbf{Example 3.6.} The hyperideal $I= 120\mathbb{Z}=\{120n:n\in\mathbb{Z}\} $  of the hyperring $ (\mathbb{Z},+,\circ) $ in  Example 3.3\textbf{(2) }is  2-absorbing primary hyperideal, but  $ \sqrt{I}=15\mathbb{Z}  $ is not prime hyperideal of $ \mathbb{Z} $. Clearly  $ 3\circ5=\{30,60\}\subseteq15\mathbb{Z} $ but $ 3,5\notin15\mathbb{Z} $ . Thus $ 15\mathbb{Z}  $ is not a prime hyperideal of $ \mathbb{Z} $.\\
\\
	\textbf{Theorem 3.8.} Let  $ P $ be a hyperideal of  $ R $ and $ I_{1},I_{2},...,I_{n} $ be 2-absorbing primary hyperideals of $ R $ such that  $ \sqrt{I_{i}}=P $ for all i=1,...,n. Then     $ \bigcap_{i=1}^{n}I_{i} $ is  2-absorbing primary hyperideal and  $ \sqrt{\bigcap_{i=1}^{n}I_{i}}=P $.\\
	\\
	\textbf{Proof.}  Clearly,  $ \sqrt{I} =\sqrt{\bigcap_{i=1}^{n}I_{i}}=\bigcap_{i=1}^{n}\sqrt{I_{i}}=P $. Let  $ I=\bigcap_{i=1}^{n}I_{i} $. Suppose that   $ abc\subseteq I $  and $ ab\nsubseteq I $ for any $ a,b,c\in I $. Hence  $ ab \nsubseteq I_{i} $ for some $ i $. Since $ I_{i}  $ is a  2-absorbing primary hyperideal  and $  abc\subseteq I\subseteq I_{i} $, then $ ac\subseteq   \sqrt{I_{i}}=P $ or $ bc\subseteq  \sqrt{I_{i}}=P $. Thus we conclude $ ac\subseteq \sqrt{I} $ or $ bc\subseteq \sqrt{I} $. Thus $ I $ is a 2-absorbing  primary hyperideal of $ R $ .$ \Box $\\
	\\
	Note that if  $  I $ and $ J $ are 2-absorbing primary hyperideals of $ R $ and    $ \sqrt{I}\neq \sqrt{J} $, then  $ I\cap J $ may not be a  2-absorbing primary hyperideal of $ R $.   We have the following example.\\
	\\
	\textbf{Example 3.9.} Consider the hyperring $ (\mathbb{Z},+,\circ) $ in Example 3.3\textbf{(1)}.  $ I=12\mathbb{Z}=\{12n:n\in\mathbb{Z}\} $ and $  J=20\mathbb{Z}=\{20n:n\in\mathbb{Z}\} $ are  2-absorbing primary hyperideal of $ \mathbb{Z} $. Since  $ \sqrt{I}=6\mathbb{Z}=\{6n:n\in\mathbb{Z}\}   $ and $ \sqrt{J}=10\mathbb{Z}=\{10n:n\in\mathbb{Z}\}  $, then  $  I\cap J=30\mathbb{Z}=\{30n:n\in\mathbb{Z}\}  $ and $ 30\mathbb{Z} $ is not a 2-absorbing primary hyperideal.\\
	\\
	\textbf{Lemma 3.10.} If $ P_{1} $ and  $ P_{2} $ are prime hyperideals of  $ R $, then $ P_{1}\cap P_{2} $ is  2-absorbing hyperideal of $ R $.\\
	\\
	\textbf{Proof.} Let $ a,b,c\in R $ such that  $ abc\subseteq  P_{1}\cap P_{2} $, $ ab\nsubseteq P_{1}\cap P_{2}  $ and $ bc\nsubseteq P_{1}\cap P_{2} $. Then $ a,b,c\notin P_{1}\cap P_{2} $. Asume that  $ a\in P_{1}\cap P_{2} $, then $ a\in P_{1} $ and $ a\in P_{2} $. Since  $ P_{1} $ and  $ P_{2} $ are hyperideals, we have  $ ab\subseteq P_{1} $ and $ ab\subseteq  P_{2} $. Then $ ab\subseteq P_{1}\cap P_{2} $, which is a contradiction. Thus  $ a\nsubseteq P_{1}\cap P_{2} $. Similarly, $ b,c\notin  P_{1}\cap P_{2} $. We consider three cases. \\
	\textbf{Case one:} Suppose that  $ a\notin P_{1} $ and $ a\notin P_{2} $. Since  $c\notin P_{1}\cap P_{2} $, we have three cases again. Assume that  $ c\notin P_{1} $ and $ c\notin P_{2} $. Since   $ P_{1} $ is a prime hyperideal and  $ ac\nsubseteq P_{1} $, $ abc\subseteq P_{1} $, then   $ b\in	P_{1} $ by Lemma 2.1. Hence $ ab\subseteq P_{1} $. Similarly,  since  $ P_{2} $  is a prime hyperideal and $ ac\nsubseteq P_{2} $, $ abc\subseteq P_{2} $, we have  $ b\in	P_{2} $ by Lemma 2.1. Hence  $ ab\subseteq P_{2} $. Thus   $ ab\subseteq P_{1}\cap P_{2} $ which is a contradiction. Thus $ c\in P_{1} $ or $ c\in P_{2} $.  Now, assume that $ c\notin P_{1} $ and $ c\in P_{2} $. Since  $ P_{1} $  is a prime hyperideal and  $ ac\nsubseteq P_{1} $, $ abc\subseteq P_{1} $ , we have   $ b\in	P_{1} $. Thus $ bc\subseteq P_{1} $. Since $ c\in P_{2} $, then   $ bc\subseteq P_{2} $ and so $bc\subseteq P_{1} \cap P_{2}$, a contradiction. Finally, assume that $ c\notin P_{2} $ and $ c\in P_{1} $. Since  $ P_{2} $  is a prime ideal and $ ac\nsubseteq P_{2} $, $ abc\subseteq P_{2} $, then $ b\in	P_{2} $ and so $ bc\subseteq P_{2} $. Since $ c\in P_{1} $, we conclude   $ bc\subseteq P_{1} $. So we have $bc\subseteq P_{1} \cap P_{2}$ which is a contradiction. Thus, if  $ a\in P_{1}\cap P_{2} $, implies that  $ a\in P_{1}$ or $  a \in P_{2} $.\\
	\textbf{Case Two:} Suppose that $ a\in P_{1} $ and $ a\notin P_{2} $. We show that $ c\in P_{2} $. Assume that $ c\notin P_{2} $. Since $ P_{2} $ is a prime hyperideal, we have $ ac\nsubseteq  P_{2} $. Since whenever $ abc\subseteq  P_{2}  $,  $ ac\nsubseteq  P_{2} $  and also  $ P_{2} $ is a prime hyperideal, then  $ b\in P_{2} $ by Lemma 2.1. Hence $ ab\subseteq P_{1}\cap P_{2} $ which is a contradiction. Thus $ c\in P_{2} $. Since   $ c\notin P_{1}\cap P_{2} $, we get  $ c\notin P_{1} $. Therefore, $ ac\subseteq P_{1}\cap P_{2} $. \\
	\textbf{Case Three:} Suppose that $ a\in P_{2} $ and $ a\notin P_{1} $. We show that $ c\in P_{1} $. Assume that $ c\notin P_{1} $. Since $ P_{2} $ is a prime hyperideal, then $ ac\nsubseteq  P_{1} $. Since whenever $ abc\subseteq  P_{1}  $, $ ac\nsubseteq  P_{1} $ and  $ P_{1} $ is a prime hyperideal, by Lemma 2.1  $ b\in P_{1} $. Hence $ ab\subseteq P_{1}\cap P_{2} $ which is a contradiction.  Since $ c\in P_{1} $ and $ c\notin P_{1}\cap P_{2} $, we have $ c\notin P_{2} $ and hence $ ac\subseteq P_{1}\cap P_{2} $. Consequently, $ P_{1}\cap P_{2} $ is 2-absorbing hyperideal. $ \Box $ \\
	\\
	The following example shows that the converse of this lemma is not true in general.\\
	\\
	\textbf{Example 3.11.} In  Example 3.3\textbf{(3)},     $ \{\bar{0}\}\cap \{\bar{0},\bar{2},\bar{4}\}=\{\bar{0}\} $. Clearly  $ \{\bar{0}\} $ is a 2-absorbing hyperideal of $ (\mathbb {Z}_6,\oplus,\star)  $, but  it is not a prime hyperideal.\\
	\\
	\textbf{Theorem 3.12.} Let $ I $ be $ P_{1} $-primary $ C $-ideal and  $ J $ be $ P_{2} $-primary  $ C $-ideal of  $ R $. Then the following statements hold.
	\begin{enumerate}
		\item $ I_{1}\cap I_{2} $ is 2-absorbing primary hyperideal.
		\item $ I_{1}I_{2} $ is 2-absorbing primary hyperideal.
	\end{enumerate} 
	\textbf{Proof.}\\
	1. Let  $ I_{1}\cap I_{2}=K $. Then $ \sqrt{K}=P_{1}\cap P_{2} $. Now, we show that  $K  $ is a 2-absorbing primary hyperideal of $ R $. Suppose that   $ abc\subseteq K $, $ ac\nsubseteq\sqrt{K} $ and $ bc\nsubseteq\sqrt{K} $ for any $ a,b,c\in R $. Since $ \sqrt{K} $ is a hyperideal, we have $ a\notin\sqrt{K} $,  $ b\notin\sqrt{K} $ and $ c\notin\sqrt{K} $. Since  $ I_{1} $ and $ I_{2} $ are  $ C $-ideals, we have $ P_{1} $ and $ P_{2} $ are prime  hyperideals. By Lemma 3.10, we conclude that $ P_{1}\cap P_{2} $ is a 2-absorbing hyperideal of $ R $. Since $  \sqrt{K}=P_{1}\cap P_{2} $ is a 2-absorbing hyperideal, then  $ ab\subseteq\sqrt{K}\subseteq P_{1} $. Since $ P_{1} $ is a prime hyperideal, $ a\in P_{1} $ or $ b\in P_{1} $. We may assume that $ a\in P_{1} $. Hence  $ a\notin P_{2} $ since $ a\notin\sqrt{K} =P_{1}\cap P_{2}$. Similarly, one can easily show that   $ b\notin P_{1} $. We claim that  $ a\in  I_{1}  $ and 
	$ b\in  I_{2}  $.  Suppose that $ a\notin I_{1} $. Since $ I_{1} $ is  $ P_{1} $-primary hyperideal and $  a\notin I_{1} $, we get  $ bc\in P_{1} $ and so  $ bc\subseteq P_{1} $. Thus  $ bc\subseteq\sqrt{K} $ and this is a contradiciton. Hence $ a\in I_{1} $. Similarly, let $ b\notin I_{2} $. Since  $ I_{2} $ is a $ P_{2} $-primary hyperideal and $ b\notin I_{2} $, we conclude that   $ ac\subseteq P_{2} $  by Lemma 2.2.  Hence  $ ac\subseteq \sqrt{K} $ since $ ac\subseteq P_{2} $ and $ a\in P_{1} $, a contradiction. Thus $ b\in I_{2} $. Therefore,  $ ab\in I_{1}\cap I_{2}= K $ .\\
	2. For any $ a,b,c\in R $, let  $ abc\subseteq I_{1}I_{2} $ and $ ab,bc\nsubseteq \sqrt{I_{1}I_{2}}=P_{1}\cap P_{2}$. Then $ a,b,c\notin   \sqrt{I_{1}I_{2}}=P_{1}\cap P_{2} $. Moreover,  we have $ ac\subseteq\sqrt{I_{1}I_{2}}=P_{1}\cap P_{2} $ since  $ P_{1}\cap P_{2}$ is  2-absorbing hyperideal. We show that $ ac\subseteq I_{1}I_{2} $. Since $ ac\subseteq \sqrt{I_{1}I_{2}}=P_{1}\cap P_{2}\subseteq P_{1} $ and $ P_{1} $ is a prime hyperideal, We get  $ a\in P_{1} $ or $ c\in P_{1} $. We may assume that  $ a\in P_{1} $. Since   $ a\notin P_{1}\cap P_{2} $, we have $ a\notin P_{2} $. Also  $ c\in P_{2} $ and  $ c\notin P_{1} $ since $ P_{2} $ is a prime hyperideal and  $ ac\subseteq \sqrt{I_{1}I_{2}}=P_{1}\cap P_{2}\subseteq P_{2} $. Now, we claim that  $ a\in I_{1} $ and $ c\in I_{2} $. Suppose that $ a\notin I_{1} $. Since $ I_{1} $ is primary hyperideal, whenever $ abc\subseteq I_{1} $ and $a \notin I_{1} $,then $ bc\subseteq\sqrt{I_{1}}=P_{1} $. Since $ c\in P_{2} $, we have $ bc\subseteq P_{1}\cap P_{2}=\sqrt{I_{1}I_{2}} $ which  is a contradiction. Thus $ a\in I_{1} $. Similarly, we conclude that $ c\in I_{2} $. Consequently, we get $ ac\subseteq I_{1}I_{2} $.$ \Box $\\
	\\
	The above theorem is not true for a finite number of hyperideals in general.\\
	\\
	\textbf{Example 3.13.} The hyperideals  $ 3\mathbb{Z}=\{3n:n\in\mathbb{Z}\} $,  $ 5\mathbb{Z}=\{5n:n\in\mathbb{Z}\} $ and  $ 7\mathbb{Z}=\{7n:n\in\mathbb{Z}\} $  in Example 3.3\textbf{(2)} are prime $ C $-ideals and so primary  $ C $-ideals such that $ \sqrt{3\mathbb{Z}}=3\mathbb{Z} $, $ \sqrt{5\mathbb{Z}}=5\mathbb{Z} $ and $ \sqrt{7\mathbb{Z}}=7\mathbb{Z} $, but $ 3\mathbb{Z}\cap 5\mathbb{Z}\cap 7\mathbb{Z}=105\mathbb{Z}  $ is not a  2-absorbing primary hyperideal.\\
	\\
	\textbf{Theorem 3.14.} If $ f:R\rightarrow S $ is a good homomorphism of hyperrings and  $ J $ is a 2-absorbing primary hyperideal of  $ S $, then   $ f^{-1}(J) $ is a 2-absorbing primary hyperideal of  $ R $.\\
	\\
	\textbf{Proof.} Let $ xyz\subseteq f^{-1}(J) $ for any $ x,y,z\in R $. Since $ f(xyz)=f(x)f(y)f(x)\subseteq J $ and $ J $ is a 2-absorbing primary  hyperideal, we have   $ f(x)f(y)\subseteq J $ or $ f(x)f(z)\subseteq \sqrt{J} $ or $ f(y)f(z)\subseteq \sqrt{J} $. Hence 
	$f(xy)\subseteq J $ or $ f(xz)\subseteq \sqrt{J} $ or $ f(yz)\subseteq \sqrt{J} $. Thus $ xy\subseteq f^{-1}(J) $ or $ xz\subseteq f^{-1}(\sqrt{J}) $ or $ yz\subseteq f^{-1}(\sqrt{J}) $. By the equality $\sqrt{f^{-1}(J)}=f^{-1}(\sqrt{J})  $, $f^{-1}(J) $ is a  2-absorbing primary  hyperideal.$ \Box $\\
	\\ 
	\textbf{ Theorem 3.15.} Let $ f:R\rightarrow S $ be an onto good homomorphism of hyperrings and $ I $ be a $ C_{u} $-ideal of the hyperring $ R $. If  $ I $ is a   2-absorbing primary hyperideal of $ R $ such that $ Kerf\subseteq I $, then   $ f(I) $ is a 2-absorbing primary hyperideal of  $ S $.\\
	\\
	\textbf{Proof.} Let  $ s_{1}s_{2}s_{3}\subseteq f(I ) $ for any  $ s_{1},s_{2},s_{3}\in S $. Since  $ f $ is onto, there exist  $ r_{1},r_{2},r_{3}\in R $ such that $s_{1}=f(r_{1}) $ and $ s_{2}=f(r_{2}) $ and $ s_{3}=f(r_{3}) $. Then $  f(r_{1}r_{2}r_{3})= f(r_{1})f(r_{2})f(r_{3})=s_{1}s_{2}s_{3}\subseteq f(I)  $ and hence $ 0\in  f(I)-f(r_{1}r_{2}r_{3})=f(I-r_{1}r_{2}r_{3})=\{f(u):u\in I-r_{1}r_{2}r_{3}\} $. Thus    $ f(v)=0\in <0> $ for some  $ v\in I-r_{1}r_{2}r_{3} $. Hence  $ v\in Kerf\subseteq I $ and so  $ I\cap (I-r_{1}r_{2}r_{3})\neq\varnothing $. Since $ I $ is a  $ C_{u} $-ideal, we have $ I-r_{1}r_{2}r_{3}\subseteq I $. Thus $ r_{1}r_{2}r_{3}\subseteq I $. Since  $ I $ is a  2-absorbing primary hyperideal, we get $ r_{1}r_{2}\subseteq I  $ or $ r_{1}r_{3} \subseteq \sqrt{I} $ or $ r_{2}r_{3}\subseteq \sqrt{I} $. By Proposition 2.8, $ s_{1}s_{2}\subseteq f(I) $ or $ s_{1}s_{3}\subseteq f(\sqrt{I})\subseteq\sqrt{f(I)} $ or $ s_{2}s_{3}\subseteq f(\sqrt{I })\subseteq\sqrt{f(I)}   $. $ \Box$\\
	\\
	\textbf{Corollary 3.16.} Let $ I $ be a   $ C_{u} $-ideal of the hyperring $ R $  and $ J $ be a hyperideal of  $ R $ such that $ J\subseteq I $. If $ I $ is a  2-absorbing primary hyperideal of $ R $, then   $ I\diagup J $ is a  2-absorbing primary hyperideal of $ R\diagup J $.\\
	\\
	\textbf{Proof.} A mapping  $ f:R\longrightarrow R/J $ with  $ f(x)=x+J $ for all $ x\in R $ is an onto good homomorphism. Since the zero in   $ R/J $ is  $ J $ and  $ <J>=\{J\} $ [see 4] and also  $ Kerf=f^{-1}(<J>)=f^{-1}(J)=\{r\in R:f(r)=J\}=\{r\in R:r+J=J\} $, then $ Kerf=J  $. Hence,  the proof is completed by  Theorem 3.15. $ \Box $\\
	\\
		\textbf{Lemma 3.17.} Let $ I $ be a 2-absorbing primary hyperideal of a strongly distributive multiplicative hyperring $ R $ and $ J $ be a hyperideal of  $ R $. If   $ abJ\subseteq I $ and $ ab\nsubseteq I $ for any $ a,b\in R $, then  $ aJ\subseteq \sqrt{I} $ or $ bJ\subseteq \sqrt{I} $.\\
		\\
		\textbf{Proof.} Suppose  that   $ aJ\nsubseteq \sqrt{I} $ ve $ bJ\nsubseteq \sqrt{I} $ for some $ a,b\in R $. Since $ aJ=\cup aj_{i}\nsubseteq \sqrt{I} $ and  $ bJ=\cup bj_{i}\nsubseteq \sqrt{I} $ for all $ j_{i}\in J $, there exists a  $ j_{i} $ such that $ aj_{i}\nsubseteq \sqrt{I} $ and $ bj_{i}\nsubseteq \sqrt{I} $. 	We may assume that $ aj_{1}\nsubseteq \sqrt{I} $ and $ bj_{2}\nsubseteq \sqrt{I} $ for some $ j_{1},j_{2}\in J $. Also for all $ j_{i} $,  $ abj_{i} \subseteq I $. Since $ abj_{1}\subseteq I $,  $ ab\nsubseteq I   $ and $ aj_{1}\nsubseteq \sqrt{I} $, we have  $ bj_{1}\subseteq \sqrt{I} $. Similarly, since $ abj_{2}\subseteq I $, $ ab\nsubseteq I   $ and $ bj_{2}\nsubseteq \sqrt{I} $, we get $ aj_{2}\subseteq \sqrt{I} $. Now, since  $ I $ is a   2-aborbing primary hyperideal, whenever 	$ ab(j_{1}+j_{2})\subseteq I $ and  $ ab\nsubseteq I $, we have  $ a(j_{1}+j_{2}) \subseteq \sqrt{I} $ or $ b(j_{1}+j_{2})\subseteq \sqrt{I} $. Suppose that $ a(j_{1}+j_{2})=aj_{1}+aj_{2}\subseteq \sqrt{I} $. Since $  aj_{2}\subseteq \sqrt{I} $, we conclude $ aj_{1}\subseteq \sqrt{I} $ which is a contradiction.  Similarly, let $ b(j_{1}+j_{2})=bj_{1}+bj_{2}\subseteq \sqrt{I} $. Since $  bj_{1}\subseteq \sqrt{I} $, we have  $ bj_{2}\subseteq \sqrt{I} $ which is a contradiction again. Thus  $ aJ\subseteq \sqrt{I} $ or $ bJ\subseteq \sqrt{I} $.$ \Box $\\
		\\
		\textbf{Theorem 3.18.} Let $ I $ be a  hyperideal of a strongly distributive multiplicative hyperring  $ R $. Then  $ I $ is a  2-aborbing primary hyperideal of   $ R $ if and only if   $ I_{1}I_{2}I_{3}\subseteq I $, then   $ I_{1}I_{2}\subseteq I $ or $ I_{2}I_{3}\subseteq \sqrt{I} $ or $ I_{1}I_{3}\subseteq \sqrt{I} $ for any  hyperideals $ I_{1},I_{2},I_{3} $ .\\
		\\
		\textbf{Proof.}  Let $ I $ be a  2-aborbing primary hyperideal and  $ I_{1}I_{2}I_{3}\subseteq I $ and  $ I_{1}I_{2}\nsubseteq I $. We claim that $ I_{2}I_{3}\subseteq \sqrt{I} $ or $ I_{1}I_{3}\subseteq \sqrt{I} $. Suppose that $ I_{2}I_{3}\nsubseteq \sqrt{I} $ and $ I_{1}I_{3}\nsubseteq \sqrt{I} $. Hence  $ i_{1}I_{3}\nsubseteq \sqrt{I} $ and $ i_{2}I_{3}\nsubseteq \sqrt{I} $  for some $ i_{1}\in	I_{1} $ and $ i_{2}\in	I_{2} $. By  Lemma 3.17, we get  $ i_{1}i_{2}\subseteq I $. Since  $ I_{1}I_{2}\nsubseteq I $,   $ ab\nsubseteq I $ for some  $ a\in I_{1} $ and $ b\in I_{2} $. Since $ abI_{3}\subseteq I_{1}I_{2}I_{3}\subseteq I $ and $ ab\nsubseteq I  $, by  Lemma 3.17  $ aI_{3}\subseteq\sqrt{I} $ or $ bI_{3}\subseteq\sqrt{I} $.\\
		\\
		\textbf{Case one:} Suppose that $ aI_{3}\subseteq\sqrt{I} $ and $ bI_{3}\nsubseteq\sqrt{I} $. Since $ i_{1}bI_{3}\subseteq I_{1}I_{2}I_{3}\subseteq I $, $ bI_{3}\nsubseteq\sqrt{I} $ and $ i_{1}I_{3}\nsubseteq\sqrt{I} $, we have $ i_{1}b\subseteq I $ by Lemma 3.17. Since $ (a+i_{1})bI_{3}\subseteq I_{1}I_{2}I_{3}\subseteq I $ and $ bI_{3}\nsubseteq\sqrt{I} $, we have  $ (a+i_{1})I_{3}\subseteq \sqrt{I} $ or $ (a+i_{1})b\subseteq I $ by Lemma 3.17. Suppose that $ (a+i_{1})I_{3}\subseteq \sqrt{I} $. Then for every $ i_{3}\in I_{3} $, since  $ R $ is  strongly distributive, we conclude that   $ (a+i_{1})I_{3}=\cup (a+i_{1})i_{3}=ai_{3}+i_{1}i_{3}=aI_{3}+i_{1}I_{3}\subseteq\sqrt{I} $. Since  $ \sqrt{I} $ is a hyperideal and $ aI_{3}\subseteq\sqrt{I} $, we get  $ i_{1}I_{3}\subseteq\sqrt{I} $ which is a contradiction. Now, suppose that  $ (a+i_{1})b=ab+i_{1}b\subseteq I $. Since $ I $ is a hyperideal and $ i_{1}b\subseteq I $, we have  $ ab\subseteq I $ which is a contradiciton again.\\
		\textbf{Case Two:} Suppose that $ aI _{3}\nsubseteq\sqrt{I}  $ and $ bI _{3}\subseteq\sqrt{I} $. Then  $ ai_{2}\subseteq I $ by Lemma 3.17. Since $ a(b+i_{2})I_{3}\subseteq I_{1}I_{2}I_{3}\subseteq I $ but $ aI_{3}\nsubseteq\sqrt{I} $,  we have   $ a(b+i_{2})\subseteq I $ or $ (b+i_{2})I_{3}\subseteq \sqrt{I} $ by Lemma 3.17. Suppose that $ (b+i_{2})I_{3}\subseteq \sqrt{I} $. Since $ R $ is strongly distributive, we get    $ (b+i_{2})I_{3}=\cup (b+i_{2})i_{3}=bi_{3}+i_{2}i_{3}=bI_{3}+i_{2}I_{3}\subseteq\sqrt{I} $ for every $ i_{3}\in I_{3} $. Since $ \sqrt{I} $ is a hyperideal  and $ bI_{3}\subseteq\sqrt{I} $, we conclude  $ i_{2}I_{3}\subseteq\sqrt{I} $ which is a contradiction. Now,  suppose that $ a(b+i_{2})b=ab+i_{2}b\subseteq I $. Similarly, since $ I $ is a hyperideal and $ ai_{2}\subseteq I $, we have  $ ab\subseteq I $ which is a contradiction.\\
		\textbf{Case Three:} Suppose that $ aI _{3}\subseteq\sqrt{I}  $ and $ bI _{3}\subseteq\sqrt{I} $. Since  $ bI_{3}\subseteq\sqrt{I} $ and $ i_{2}I_{3}\nsubseteq\sqrt{I} $, we have $ (b+i_{2})I_{3}\nsubseteq\sqrt{I} $. By Lemma 3.17, we conclude that  $ i_{1}(b+i_{2})=i_{1}b+i_{1}i_{2}\subseteq I $. Since  $ i_{1}i_{2}\subseteq I $ and $ i_{1}b+i_{1}i_{2}\subseteq I $, we get $ bi_{1}\subseteq I $. Since  $ aI_{3}\subseteq\sqrt{I} $ and $ i_{1}I_{3}\nsubseteq\sqrt{I} $, we conclude that $ (a+i_{1})I_{3}\nsubseteq\sqrt{I} $.  Hence   $ (a+i_{1})i_{2}=ai_{2}+i_{1}i_{2}\subseteq I $ by Lemma Teorem 3.17. Since  $ i_{1}i_{2}\subseteq I $ and $ ai_{2}+i_{1}i_{2}\subseteq I $, we have  $ ai_{2}\subseteq I $. Thus    $ (a+i_{1})(b+i_{2})=ab+ai_{2}+bi_{1}+i_{1}i_{2}\subseteq I $ by  Lemma 3. 17. Since    $ ai_{2},bi_{1},i_{1}i_{2}\subseteq I $,  we have $ ai_{2}+bi_{1}+i_{1}i_{2}\subseteq I $. Hence $ ab\subseteq I $ since  $ ab+ai_{2}+bi_{1}+i_{1}i_{2}\subseteq I $ and $ ai_{2}+bi_{1}+i_{1}i_{2}\subseteq I $,a contradiction. Consequently, we conclude that $ I_{1}I_{3}\subseteq\sqrt{I} $ or $ I_{2}I_{3}\subseteq\sqrt{I} $ .\\
		\\
	Conversely, suppose that whenever $ I_{1}I_{2}I_{3}\subseteq I $, then   $ I_{1}I_{2}\subseteq I $ or $ I_{2}I_{3}\subseteq \sqrt{I} $ or $ I_{1}I_{3}\subseteq \sqrt{I} $ for any  hyperideals $ I_{1},I_{2},I_{3} $. Cleary, $ I $ is 2-absorbing primary hyperideal of $ R $ by the definition of  principal hyperideal of $ R $.	$ \Box $

\end{document}